\documentclass[12pt]{article}
\usepackage[english]{babel}
\usepackage[cp1251]{inputenc}
\usepackage{amscd}
\usepackage{hyperref}
\usepackage{amsmath}
\usepackage{amssymb}

\newtheorem{thm}{Theorem}
\newtheorem{lem}{Lemma}
\newtheorem{propos}{Proposition}
\newtheorem{cor}{Corollary}
\begin{document}

%\begin{frontmatter}

\noindent{\Large\bf Rate of approximation of $zf'(z)$ by special sums
associated with the zeros of the Bessel polynomials}

\bigskip

\noindent{\bf Mikhail~A.~Komarov} (\verb"kami9@yandex.ru")

\noindent{Vladimir State University, 600000 Vladimir, Russia}

\bigskip

\noindent{\bf Abstract.}
Let $\alpha_{n1},\dots,\alpha_{nn}$ be the zeros of the $n$th
Bessel polynomial $y_n(z)$ and let $a_{nk}=1-\alpha_{nk}/2$,
$b_{nk}=1+\alpha_{nk}/2$ $(k=1,\dots,n)$. We propose the new
formula \[z f'(z)\approx \sum_{k=1}^n \big(f(a_{nk} z)-f(b_{nk}
z)\big)\] for numerical differentiation of analytic functions
$f(z)=\sum_0^\infty f_m z^m$. This formula is exact for all
polynomials of degree at most $2n$. We find the sharp order of
nonlocal estimate of the corresponding remainder for the case when
all $|f_m|\le 1$. The estimate shows a high rate of convergence of
the differentiating sums to $zf'(z)$ on compact subsets of the
open unit disk, namely, $O(0.85^n n^{1-n})$ as $n\to \infty$.



%\end{frontmatter}

%\linenumbers

\section{Introduction}\label{Sec 1}

Throughout the paper we assume that $f$ is a function, analytic in
a neighborhood of the origin, and put $f(z)=\sum_0^\infty f_m
z^m$. Consider the elementary 2-node differentiation formula
\begin{equation}\label{a=3/2, b=1/2}
    zf'(z)\approx f(3z/2)-f(z/2).
\end{equation}
We can easily seen that
\[zf'(z)-\{f(3z/2)-f(z/2)\}=O(z^3) \qquad {\rm as} \quad z\to 0,\]
whereas for any complex numbers $a, b$ such that $a\ne 3/2$ and/or
$b\ne 1/2$ the difference $zf'(z)-\{f(az)-f(bz)\}$ is only
$O(z^\gamma)$, $\gamma\le 2$. Hence the formula (\ref{a=3/2,
b=1/2}) is exact for all polynomials of degree at most $2$ and is
optimal (by the order of local approximation) in the class of all
formulas of the form $zf'(z)\approx f(az)-f(bz)$. Moreover, the
nodes $3z/2$ and $z/2$ in (\ref{a=3/2, b=1/2}) are independent of
$f$.

In this note we construct the $2n$-node generalization of
(\ref{a=3/2, b=1/2}) which preserves the main features of
(\ref{a=3/2, b=1/2}) --- the exactness for polynomials of degree
at most $2n$ ($2n$ is the number of nodes), the property that the
nodes are independent of $f$ and that the values of $f$ in them
are summed with the multipliers $+1$ or $-1$ (some known similar
differentiation formulas are discussed in Sec. \ref{Sec 4}).

We also find the sharp in order of the quantity $n$ nonlocal
estimate of the corresponding remainder for the case of bounded
Taylor coefficients, $f_m$. The estimate shows a very high rate of
convergence of the differentiating sums to $zf'(z)$ on compact
subsets of the disk $|z|<1$, namely, $O(0.85^n n^{1-n})$ as $n\to
\infty$.

The nodes in our differentiation formula linearly depend on the
zeros $\alpha_{nk}$ ($k=1,\dots,n$) of the $n$th Bessel
polynomial $y_n(z)$, first introduced in \cite{KrallFrink}:
\[y_n(z)=\sum_{k=0}^{n} \frac{(n+k)!}{(n-k)!k!}\left(\frac{z}{2}\right)^k.\]
Recall that all $\alpha_{nk}$ are simple \cite[p.\,75]{Grosswald}
and lie in the semi-annulus
\begin{equation}\label{semi-ann}
    \frac{1}{n+2/3}\le |\alpha_{nk}|\le \frac{2}{n+1}, \qquad
    {\rm Re}\,\alpha_{nk}<0
\end{equation}
\cite[Sec. 5]{BSV}, and there are algorithms for the accurate
computation of these zeros \cite{Pasquini}. We also need the
following well-known results on the power sums of $\alpha_{nk}$:
if $j=1,2,\dots$ and
$\sigma^n_j=\alpha_{n1}^j+\dots+\alpha_{nn}^j$, then
\begin{align}
 &   \sigma_1^{n}=-1, \qquad \sigma_{2j+1}^{n}=0 \quad
    (j=1,\dots,n-1); \label{sigm_1=-1...}\\
 &   \sigma^{n}_{2n+1}=(-1)^n 4^n n!^2/(2n)!^2. \label{sigm 2n+1}
\end{align}
The formulas (\ref{sigm_1=-1...}) and (\ref{sigm 2n+1}) were
proved in \cite{Burch} and in \cite{IsmailKelker}, respectively.

\section{The differentiation formula and its local properties}\label{Sec 2}

For $n=1,2,\dots$, define the quantities
\[a_{nk}=1-\alpha_{nk}/2, \qquad b_{nk}=1+\alpha_{nk}/2 \qquad (k=1,\dots,n)\]
and introduce the main differentiation formula:
\begin{equation}\label{zf' approx sum}
    z f'(z)\approx \sum_{k=1}^n \big(f(a_{nk} z)-f(b_{nk} z)\big).
\end{equation}
Since $y_1(z)=1+z$, then $\alpha_{1,1}=-1$, such that the formula
(\ref{a=3/2, b=1/2}) is actually the particular case of (\ref{zf'
approx sum}) with $n=1$.

Find the local representation of the remainder
\[r_n(f;z):=z f'(z) - \sum_{k=1}^n \big(f(a_{nk} z)-f(b_{nk} z)\big)\]
of (\ref{zf' approx sum}) near the origin. Obviously, we have
\begin{align}
r_n(f;z)&=\sum_{m=1}^\infty m f_m z^m-\sum_{k=1}^n
\sum_{m=1}^\infty (a_{nk}^m-b_{nk}^m) f_m z^m \notag\\
&=\sum_{m=1}^\infty (m-A_m^{n}) f_m z^m, \qquad {\rm where}\qquad
A_m^{n}:=\sum_{k=1}^n \big(a_{nk}^m -
b_{nk}^m\big).\label{rn(f;z)}
\end{align}

\begin{propos}\label{Prop 1}
For $n=1,2,\dots$,
\[m-A_m^n=0 \quad (m=1,\dots,2n); \qquad (2n+1)-A_{2n+1}^n=
(-1)^n \frac{n!^2}{(2n)!^2}.\]
\end{propos}

{\sc Proof.} For $m=1,2,\dots$, by applying the binomial formula, we have
\[(1-t)^m-(1+t)^m=-2\sum_{j=0}^{M} \binom{m}{2j+1}t^{2j+1},\]
where $M=[(m-1)/2]$ and $[x]$ denotes the greatest integer $\le
x$, so that
\[A_m^{n}=\sum_{k=1}^n \left\{(1-\alpha_{nk}/2)^m-(1+\alpha_{nk}/2)^m\right\}=
-\sum_{j=0}^{M} \binom{m}{2j+1} \frac{\sigma^{n}_{2j+1}}{4^j}.\]
If $m\le 2n$, then $M\le n-1$ and $A_m^{n}=-\binom{m}{1}
\sigma^{n}_{1}-0=m$ by (\ref{sigm_1=-1...}). Analogously, if
$m=2n+1$, then $M=n$, so that
\[A_{2n+1}^{n}=-\binom{2n+1}{1}
\sigma^{n}_{1}-0-\binom{2n+1}{2n+1}\frac{\sigma^{n}_{2n+1}}{4^n}=
(2n+1)-(-1)^n\frac{n!^2}{(2n)!^2}\] by (\ref{sigm_1=-1...}) and
(\ref{sigm 2n+1}). This completes the proof of Proposition
\ref{Prop 1}.

\medskip

\begin{cor}
For $n=1,2,\dots$, the remainder of $(\ref{zf' approx sum})$ has
the local form
\[r_n(f;z)=(-1)^n\frac{n!^2}{(2n)!^2}\cdot
f_{2n+1}z^{2n+1}+O(z^{2n+2}) \qquad as\quad z\to 0.\] In
particular, the sum in the right side of $(\ref{zf' approx sum})$
interpolates the function $zf'(z)$ at zero with the multiplicity
$2n+1$ and $r_n(p;z)\equiv 0$ for all polynomials $p$ of degree at
most $2n$.
\end{cor}

We conclude this section with some properties of the quantities
$a_{nk}$, $b_{nk}$.

Because of the properties of $\alpha_{nk}$ (see Sec. \ref{Sec 1}),
{\it all $a_{nk}$ and $b_{nk}$ $(k=1,\dots,n)$ are pairwise
distinct and tend to $1$ as $n\to \infty$}. We now prove

\begin{propos}\label{Prop 1+}
For $n=1,2,\dots$, $|a_{nk}|>1$ and $|b_{nk}|<1$ $(k=1,\dots,n)$.
\end{propos}

{\sc Proof.} Since ${\rm Re}\,\alpha_{nk}<0$, then $|a_{nk}|\ge {\rm
Re}\,a_{nk}=1-{\rm Re}\,\alpha_{nk}/2>1$.

To prove $|b_{nk}|<1$ we need the more strong estimate \cite[Eq.
(6.2)]{BSV}:
\[{\rm Re}\,\alpha_{nk}<-[n^3(n+1)]^{-1/2}.\]  Thus we actually have
\begin{align*}
|b_{nk}|^2& =(1+{\rm Re}\,\alpha_{nk}/2)^2+({\rm
           Im}\,\alpha_{nk}/2)^2=1+{\rm Re}\,\alpha_{nk}+|\alpha_{nk}|^2/4 \\
          & <1-\frac{1}{\sqrt{n^3(n+1)}}+\frac{1}{(n+1)^2}=
           1+\frac{1}{(n+1)^2}\left[1-\left(\frac{n+1}{n}\right)^{3/2}\right]<1,
\end{align*}
and the assertion follows.

\section{The nonlocal estimates of the interpolation error in (\ref{zf' approx sum})}\label{Sec 3}

Put $n_0(x)=\max\left\{14;\,x^2(1-x)^{-2}\right\}-1$.

\begin{thm}\label{Th1}
Let $x\in (0,1)$ and all $|f_m|\le 1$. If $n\ge n_0(x)$, then
\begin{equation}\label{|Delta_n(f;z)|<...}
    |r_n(f;z)|\le \frac{|z|^{2n+1}}{x-|z|}\frac{2x}{(1-x)^{2n+2}}
    \frac{n!^2}{(2n)!^2}, \qquad |z|<x,
\end{equation}
and this estimate is sharp in order of the quantity $n$ for
$|z|\approx x$.

In particular, if $n\ge 13$ and $x=x_n=\sqrt{n}/(1+\sqrt{n})$,
then
\begin{equation}\label{|Delta_n(f;z)|<...(x_n)}
    |r_n(f;z)|\le \frac{|z|^{2n+1}}{x_n-|z|}\frac{0.92^{2n}}{n^{n-1}}, \qquad |z|<x_n.
\end{equation}
\end{thm}

{\it Remark 1.}
Under the condition $|f_m|\le 1$, the function $r_n(f;\cdot)$ is
analytic in the disk $|z|<(n+1)/(n+2)$. Indeed, $f$ is analytic in
the open unit disk $D$, whereas
\[|b_{nk}|<1, \qquad |a_{nk}|\le 1+|\alpha_{nk}|/2\le (n+2)/(n+1)\]
(use Proposition \ref{Prop 1+} and $|\alpha_{nk}|\le 2/(n+1)$),
hence all the nodes $a_{nk}z$ and $b_{nk}z$ belong to $D$ if
$|z|<(n+1)/(n+2)$. In particular, $r_n(f;\cdot)$ is analytic for
$|z|<x$, since by $n\ge n_0(x)$ we have $x\le
\sqrt{n+1}/(1+\sqrt{n+1})<(n+1)/(n+2)$. Note that, generally
speaking, $r_n(f;\cdot)$ is not analytic in the whole disk $D$ by
$|a_{nk}|>1$.

\medskip

{\it Remark 2.}
The choice $x_n=\sqrt{n}/(1+\sqrt{n})$ is admissible, since $n\ge
n_0(x_n)$ for $n\ge 13$.

{\sc Proof.} First consider the model function $g(z)=z/(1-z)=z+z^2+\dots$.
Obviously, $r_n(g;z)\equiv zR_n(z)$, where
\begin{equation}\label{Rn(z)}
    R_n(z)=\frac{1}{(z-1)^2}-\sum_{k=1}^n \frac{a_{nk}}{1-a_{nk} z}+
    \sum_{k=1}^n \frac{b_{nk}}{1-b_{nk} z}.
\end{equation}
The following assertion will be established in Sec. \ref{Sec 5},
and this proves, in particular, the sharpness of
(\ref{|Delta_n(f;z)|<...}) for $|z|\approx x$.

\begin{propos}\label{Prop 2}
Let $x\in (0,1)$. If $n\ge n_0(x)$, then
\begin{equation}\label{...<|Delta_n|<...}
    \frac{x^{2n}}{(1-x)^{2n+2}}\frac{n!^2}{(2n)!^2}\le
    \max_{|z|=x}|R_n(z)|<\frac{2
    x^{2n}}{(1-x)^{2n+2}}\frac{n!^2}{(2n)!^2}.
\end{equation}
\end{propos}

Now assume $x\in (0,1)$, $|z|<x$ and $n\ge n_0(x)$. Then we have
$|a_{nk}z|<1$, $|b_{nk}z|<1$ (see Remark 1) and
\[R_n(z)=\sum_{m=1}^{\infty} (m-A_m^{n}) z^{m-1}\equiv
\sum_{m=2n+1}^{\infty} (m-A_m^{n}) z^{m-1}\] (see Proposition
\ref{Prop 1}), therefore, by using the Cauchy inequalities
\[\big|m-A_m^{n}\big|\le x^{1-m}\max_{|z|=x}|R_n(z)|\] and the upper
estimate in (\ref{...<|Delta_n|<...}), we get

\begin{cor}
Let $x\in (0,1)$. If $n\ge n_0(x)$, then
\begin{equation}\label{T1-2}
    \big|m-A_m^{n}\big|<\frac{2x^{2n+1-m}}{(1-x)^{2n+2}}\frac{n!^2}{(2n)!^2}
    \qquad (m\ge 2n+1).
\end{equation}
\end{cor}

(The right side in (\ref{T1-2}) tends to infinity as $m\to
\infty$, and the example $A_m^1=(3^m-1)/2^m$ shows that the
quantities $\big|m-A_m^n\big|$ are indeed unbounded as $m\to
\infty$.)

Thus (\ref{|Delta_n(f;z)|<...}) immediately follows from
(\ref{rn(f;z)}), $|f_m|\le 1$, Proposition \ref{Prop 1} and
(\ref{T1-2}):
\[|r_n(f;z)|\le \sum_{m=2n+1}^\infty \big|m-A_m^{n}\big| |z|^m\le
\frac{2|z|^{2n+1}}{(1-x)^{2n+2}}\frac{n!^2}{(2n)!^2}
\sum_{m=0}^\infty \frac{|z|^m}{x^m}.\] To prove
(\ref{|Delta_n(f;z)|<...(x_n)}), put $x=x_n=\sqrt{n}/(1+\sqrt{n})$
in (\ref{|Delta_n(f;z)|<...}) and apply the Stirling formula
$n!=\sqrt{2\pi n} (n/e)^{n}\exp c_n$, $0<c_n<1/(12n)$, and also
the inequalities $x_n<1$, $(1+1/\sqrt n)^{\sqrt n}<e$,
$c_n+\sqrt{n}+n^{-1/2}<0.3n$ ($n\ge 13$) and $e^{1.3}/4<0.92$:
\[\frac{2x_n}{(1-x_n)^{2n+2}}\frac{n!^2}{(2n)!^2}<
\frac{(1+\sqrt n)^{2n+2}e^{2n+2c_n}}{4^{2n}n^{2n}}<
\frac{e^{2(n+c_n+\sqrt{n}+n^{-1/2})}}{4^{2n}n^{n-1}}<
\frac{0.92^{2n}}{n^{n-1}}.\] This completes the proof of Theorem
\ref{Th1}. Note that $0.92^2<0.85$.

\medskip

Analogously, by using the relation $\overline{\lim}_{m\to \infty}
|f_m|^{1/m}\le 1$ for an arbitrary analytic in $D$ function $f$,
we get the some more general result (see also \cite[Sec.
1]{D-2008}):

\begin{thm}\label{Th2}
Let $x\in (0,1)$ and $f$ is analytic in $D$. If $n\ge n_0(x)$,
then for any $\varepsilon\in (0,1)$
\[|r_n(f;z)|\le C(f,\varepsilon)\frac{(1-\varepsilon/2)^{2n}x^{2n+1}}{(1-x)^{2n+2}}
\frac{n!^2}{(2n)!^2}, \qquad |z|\le (1-\varepsilon)x.\]
\end{thm}

\section{Similar differentiation formulas}\label{Sec 4}

Putting $f(z)=f(0)+zh(z)$ in (\ref{zf' approx sum}), we get
\begin{equation}\label{H n(a)-H n(b)}
    (z h(z))'=\sum_{k=1}^n \big(a_{nk} h(a_{nk} z)-b_{nk} h(b_{nk}
    z)\big)+O(z^{2n}).
\end{equation}
The similar formula was obtained in \cite[Sec.\,2.4]{D-2008}:
\[(zh(z))'=H_n(h;z)+O(z^n), \qquad H_n(h;z):=\sum_{k=1}^n\lambda_kh(\lambda_kz).\]
Here $\{\lambda_k\}$ is the (unique) solution of the system
$\sum_1^n \lambda_k^m=m$ ($m=1,\dots,n$); $\sum_{1}^n
\nu_kh(\nu_kz)$ are so-called {\it $h$-sums}, so that the sum in
(\ref{H n(a)-H n(b)}) is the {\it difference of $h$-sums}. For the
case when $h$ is analytic in $D$, the rate of the uniform
convergence $H_n(h;z)\to (zh(z))'$ as
$n\to \infty$ on compact subsets of $D$ is geometric \cite{D-2008}.

More recently, the other similar formula was constructed in
\cite[Sec.\,5.4]{DCh-2016}:
\begin{align}
    z f'(z)&=\sum_{k=1}^n \mu_k f(\lambda_k z)-p f_{n-1}z^{n-1}-
    q f_{2n-1}z^{2n-1}+r_n(z),\label{AFS}\\
    r_n(z)&=\frac{6np}{(n-1)(n-2)}\cdot f_{2n}z^{2n}+O(z^{2n+1})
    \qquad (n\ge 3). \notag
\end{align}
Here $\{\mu_k,\lambda_k\}$ is the solution of the discrete moment
problem
\begin{equation}\label{sum mu_k lambda_k^m =c_m}
    \sum_{k=1}^n \mu_k \lambda_k^m =c_m, \qquad m=0,1,\dots,2n-1; \quad
    \mu_k\ne 0; \quad \lambda_k\ne \lambda_j \ \ \text{for} \ \ k\ne j,
\end{equation}
where $c_m=m$ ($m\ne n-1, 2n-1$), $c_{n-1}=n-1+p$,
$c_{2n-1}=2n-1+q$; $p,q$ are constant parameters for the unique
solvability of (\ref{sum mu_k lambda_k^m =c_m}) (they are
independent of $f$). The remainder in (\ref{AFS}) is of quite high
order, $O(z^{2n})$, and this is achieved by knowing only $n+2$
values of $f$ and its fixed derivatives (while the order
$O(z^{2n+1})$ in (\ref{zf' approx sum}) is achieved by knowing
$2n$ values of $f$). However, to use (\ref{AFS}) we need to
calculate these derivatives, $f^{(n-1)}(0)$ and $f^{(2n-1)}(0)$,
very accurate. We are not aware of any effective estimates
of the nonlocal error in (\ref{AFS}).

The optimal ``real'' $n$-node formula for calculation of the first
derivative of real functions at zero was obtained in
\cite{AshJansonJones} in the form
\[f'(0)x\approx \sum_{k=1}^n w_k f(u_k x) \qquad (x,w_k,u_k\in
\mathbb{R}; \ \ |u_k-u_j|\ge 1 \ {\rm for} \ k\ne j);\] here
$w_k,u_k$ satisfy $d_0=d_2=\dots=d_n=0$, $d_1=1$ and
$|d_{n+1}|={\rm minimum}$, where $d_m:=\sum_{k=1}^n w_k u_k^m$
(see also References in \cite{AshJansonJones,DCh-2016}).

\section{Proof of Proposition \ref{Prop 2}}\label{Sec 5}

Consider the normalized polynomial
\[Q_n(z)=\frac{n!}{(2n)!}\sum_{k=0}^{n} \frac{(2n-k)!}{(n-k)!}\frac{z^k}{k!}
\qquad (Q_n(0)=1),\] and let $z_{n1},\dots,z_{nn}$ be the zeros of
$Q_n$. By a simple observation,
\[Q_n(z)\equiv \frac{n!}{(2n)!}\cdot z^n y_n\left(\frac{2}{z}\right),\]
therefore $z_{nk}=2/\alpha_{nk}$ ($k=1,\dots,n$). We need

\begin{lem}[\cite{K-AlgAn}]
Let $\rho>0$. If $n\ge \max\{14; \rho^2\}-1$, then
\begin{equation}\label{|QnPn|>1/2}
    1/2<|Q_n(z)Q_n(-z)|<3/2, \qquad |z|\le \rho,
\end{equation}
and if $n>\max\{2; (\rho^2+4)/8\}$, then
\begin{equation}\label{1>|Qn(x)Pn(x)|>1-...}
    0<1-\rho^2(8n-4)^{-1}\le Q_n(x)Q_n(-x)\le 1, \qquad -\rho\le x\le
    \rho.
\end{equation}
\end{lem}

For the proof of (\ref{|QnPn|>1/2}) and
(\ref{1>|Qn(x)Pn(x)|>1-...}), see \cite[Sec.\,6]{K-AlgAn} and
\cite[Sec.\,7]{K-AlgAn}, respectively.

Define the polynomials \[M_n(z)=(1-z)^n
Q_n\left(\frac{z}{z-1}\right), \qquad G_n(z)=(1-z)^n
Q_n\left(\frac{z}{1-z}\right)\] and put $u_{nk}=a_{nk}^{-1}$,
$v_{nk}=b_{nk}^{-1}$ ($k=1,\dots,n$); the last numbers are
pairwise distinct, as well as differ of 1, since $\alpha_{nk}$ are
pairwise distinct and Re\,$\alpha_{nk}<0$.

Obviously, $\{u_{nk}\}$ and $\{v_{nk}\}$ are all the roots of
$M_n$ and $G_n$, respectively, for example,
$M_n(u_{nk})=(1-u_{nk})^n Q_n(z_{nk})=0$. Hence (\ref{Rn(z)}) takes
the form
\begin{align*}
R_n(z)&=\frac{1}{(z-1)^2}+\sum_{k=1}^n
        \bigg(\frac{1}{z-u_{nk}}-\frac{1}{z-v_{nk}}\bigg)=\frac{1}{(z-1)^2}+\frac{d}{dz}\left(\ln
        \frac{M_n(z)}{G_n(z)}\right)\\
      &=\frac{1}{(z-1)^2}-\frac{1}{(z-1)^2}\frac{d}{dw}\left(\ln
        \frac{Q_n(w)}{Q_n(-w)}\right),
\end{align*}
where $w=z/(z-1)$. Now we need the formula
\[1-\frac{d}{dz}\left(\ln \frac{Q_n(z)}{Q_n(-z)}\right)=
(-1)^{n}\frac{z^{2n}\gamma_n}{Q_n(z)Q_n(-z)}, \qquad
\gamma_n:=\frac{n!^2}{(2n)!^2}\] (see Sections 2, 5 in
\cite{K-AlgAn}; we can also obtain this formula from \cite[Eq.
(17)]{Burch}). By using this formula, we get the identity
\begin{equation}\label{Rnn-F}
    R_n(z)=(-1)^n \frac{z^{2n}\gamma_n}{(z-1)^{2n+2}}\frac{1}{Q_n(w)Q_n(-w)},
    \qquad w=z/(z-1).
\end{equation}
Let $x\in (0,1)$. If $|z|\le x$, then $|w|\le \rho_0$ with
$\rho_0:=x/(1-x)$, $\rho_0>0$. Note that $\max\{14;
\rho_0^2\}-1\equiv n_0(x)$, hence we have \[|Q_n(w)Q_n(-w)|>1/2\]
for $n\ge n_0(x)$ by (\ref{|QnPn|>1/2}). Thus the upper bound in
(\ref{...<|Delta_n|<...}) follows from (\ref{Rnn-F}).

Analogously, the lower bound in (\ref{...<|Delta_n|<...}) follows
from (\ref{Rnn-F}) and (\ref{1>|Qn(x)Pn(x)|>1-...}):
\[\max_{|z|=x}|R_n(z)|\ge |R_n(x)|=
\frac{x^{2n}\gamma_n}{(1-x)^{2n+2}}
\frac{1}{Q_n(-\rho_0)Q_n(\rho_0)}\ge
\frac{x^{2n}\gamma_n}{(1-x)^{2n+2}}\] (use that $n_0(x)=\max\{14; \rho^2\}-1>\max\{2; (\rho_0^2+4)/8\}$).

This completes the proof of Proposition \ref{Prop 2}.

\end{document}